\input amstex
\documentstyle {amsppt} \topmatter  \title The constant in the functional equation and derived exterior powers \endtitle
\vsize = 7in \hsize = 5.5 in
\author Stephen Lichtenbaum \endauthor \endtopmatter
\document \magnification = \magstep 1 \baselineskip = 1.5 \baselineskip
\define \Z {\Bbb Z}
\define \R {\Bbb R}
\define \Q {\Bbb Q}
\define \C {\Bbb C}

\S 0.  Introduction 

Let $X$ be a regular scheme, projective and flat over Spec $\Z$, of dimension $d$.  Let $f:X\to $ Spec $\Z$. The zeta-function $\zeta(X,s)$ of $X$ is defined to be $\prod_x 1/(1-(Nx)^{-s})$, where the product runs over the closed points $x$ of $X$, and $N(x)$ is the order of the residue field $\kappa (x)$.  this product is known to converge for $Re(s) > d$.  It is conjectured that $\zeta(X,s)$ can be extended to a function meromorphic in the plane,

It is further conjectured [B1] that there is a product of Gamma-functions $\Gamma(X,s)$ and a positive rational number $A$ associated with $X$ such that if $\phi(X,s) =  \zeta(X,s) \Gamma(X,s)A^{-s/2}$, then $\phi(X, s) = \pm \phi(X, d-s)$.

Let $\Omega$ be the sheaf of Kahler differentials $\Omega_{X/\Z}$ on $X$.  We are going to give a series of conjectures relating $A$ to the \'etale cohomology of the derived exterior powers of $\Omega$, and a proof of these conjectures for $d = 1$ and $d = 2$.  

We begin by reviewing the notion of derived exterior power.  If $\Cal C$ is any abelian category, there is an isomorphism of categories $N$ from the category of simplicial objects $\Cal S(\Cal C)$ of $\Cal C$ to the category $\Cal {CH}(\Cal C)$ of chain complexes of objects of $\Cal C$ whose degree is bounded below by $0$, and an explicit inverse isomorphism $K$. Now let $\Cal C$ be the category of coherent sheaves on $X$.   Let $\Lambda^k$ be the usual $k$-th exterior power on coherent sheaves.    

Let $E$ be any coherent sheaf on $X$.  Let $P^.$ be a finite resolution of $E$ by locally free sheaves, and let the $k$-th derived exterior power $\tilde \lambda^k(E)$ in the derived category of coherent sheaves on $X$ be  $N \Lambda^k KP^.$.  This is independent (in the derived category) of the choice of resolution. (See the Appendix, Theorem A.1).

It is easy to see that there is a map $\rho^k$ in the derived category from $\tilde \lambda^k E$ to $\Lambda ^k E$.  Since X is regular and projective over Spec $\Z$, we can embed $X$ as a locally complete intersection $i:X \to P$ in some projective space $P$ over $\Z$. If we let $I$ be the sheaf of ideals defining $X$, we have the exact sequence $0 \to I/I^2 \to i^*\Omega_{P/\Z} \to \Omega_{X/\Z}\to 0$, where $\Omega$ denotes as usual the Kahler differentials.

 The sheaf $I/I^2$ is locally free of rank $m$, say, and the sheaf $i^*\Omega_{P/\Z}$ is locally free of rank $n$, with $n - m = d-1$.  The canonical class $\omega = \omega_{X/\Z}$  may be defined as $\underline {Hom} ( \Lambda^m I/I^2, \Lambda^n i^*\Omega_{P/\Z})$, from which we obtain a map from $\Lambda^{d-1}\Omega_{X/\Z}$ to  $\omega$. 
 
  Composing this with $\rho^{d-1}$, we obtain a natural map $\psi$ from $\tilde\lambda^{d-1}\Omega$ to $\omega$, so we can compose the derived tensor product map from $\tilde\lambda^r\Omega \otimes _L\tilde \lambda^{d-1-r} \Omega$ to $\tilde \lambda^{d-1}\Omega$ with $\psi$ to obtain a map from $\tilde \lambda^r\Omega \otimes_L \tilde\lambda^{d-1-r}\Omega$ to $\omega$, and hence by adjointness, a map $\phi_{r,d}$ from $\tilde\lambda^r \Omega$ to $\underline{RHom}(\tilde \lambda^{d-1-r}\Omega,  \omega)$..By taking the cone of $\phi_{r,d}$ we obtain an object $C_{r,d}$ of the derived category, defined up to a non-canonical isomorphism.

By Serre duality, $\phi_{r,d}$ is an isomorphism at smooth points of $X$, so the  \`etale cohomology groups of $C_{r,d}$ have support in the bad fibres of $f$, hence are finite.  so we may define the Euler characteristic $\chi_{r,d}$  of $C_{r,d}$ to be the alternating product of the orders of those \`etale  cohomology groups.

In [L], I made a very general conjecture giving a formula for the value of the  leading term of the Laurent series expansion of $\zeta(X, s)$ at a rational  
 integer $r$.  Roughly speaking, the compatibility of this conjecture with the functional equation would result from knowing the following conjecture.

\proclaim {Conjecture 0.1} : $\chi_{r,d}$ is equal to $A$ if $r$ is odd and equal to $A^{-1}$ if r is even. \endproclaim

  In this paper we will give a proof of this result if $d$ is equal to $1$ or $2$.  The proof for $d =2$ relies heavily on deep results of Spencer Bloch [B]. I thank him for many extremely helpful discussions.  

\S 1.  Derived Tensor Product

Now assume that $X$ is an arithmetic surface, i. e. $d = 2$.

\proclaim {Theorem 1.1}  The Grothendieck group $K(X_{bad})$ of the category of coherent sheaves on $X$ with support contained in the bad fibers of $f$ is generated by  the sheaves $(i_P)_* \kappa(P)$ and $(i_Y)_*O_Y$, where $P$ is a closed point on a bad fiber , $\kappa(P)$ is the residue field at $P$, and $Y$ is an irreducible component of a bad fiber. \endproclaim

There are surjective maps $r_Y$ from $K(X_{bad})$ to $\Z$ which take the class of a sheaf $E$ to the length of $E_Q$ where $Q$ is the generic point of $Y$. 

\proclaim {Lemma 1.2} If an element $F$ of $K(X_{bad})$ is in the kernel  $Ker(X)$ of all the maps $r_Y$, then $F$ is in the subgroup of $K(X_{bad})$ generated by the residue fields. \endproclaim

\proclaim {Theorem 1.3} (computation of Euler characteristics of derived tensor products) If $P_1$ and $P_2$ are closed points on $X$ the Euler characeristic of $\kappa(P)_1 \otimes_L \kappa(P_2)$ is equal to $1$, whether or not $P_1 = P_2$.  The Euler characteristic of $O_Y \otimes _L \kappa(P)$ is also equal to $1$, whether or not $P$ is on $Y$.  The Euler characteristic of $O_{Y_1} \otimes _L O_{Y_2}$  is the intersection number $(Y_1.Y_2)$, whether or not $Y_1 = Y_2$.  \endproclaim

Proof If $P \neq Q$ the support of $\kappa(P) \otimes _L\kappa(Q)$ is contained in $P \cap Q = \phi.$ To compute $\kappa(P) \otimes_L \kappa(P) $ take the Koszul resolution $ 0 \to A \to A^2 \to A \to 0$ of $\kappa(P)$ and tensor it with $\kappa(P)$ to obtain $0 \to \kappa(P) \to \kappa(P)^2 \to \kappa(P) \to 0$ which clearly has Euler characteristic equal to $1$.

\proclaim {Corollary 1.4} If $F_1$ and $F_2$ are in $Ker(X)$,  $\chi(F_1 \otimes _L F_2 )= 1$. \endproclaim

\S 2. Derived Exterior Powers.

We now give  the basic lemma ([H] , Chapter II, Exercise 5.16) about derived exterior powers:

\proclaim{Lemma 2.1} Let $0 \to E_1 \to E_2 \to E_3 \to 0 $ be an exact sequence of coherent locally free sheaves on $X$. Let $r \geq 1$ be an integer.   Then there is a filtration on $\Lambda^rE_2$:

$$ 0 =G_0 \subseteq G_1 \subseteq \dots \subseteq G_{r-1} \subseteq G_r = \Lambda^rE_2$$

and exact sequences $ 0 \to G_{i-1} \to G_i \to \Lambda^{r-i} E_1 \otimes \Lambda^i E_3 \to 0 $ .  The filtration and the associated maps are functorial in exact   sequences.  \endproclaim

If $E$ is a coherent sheaf on $X$, let $[E]$ denote the class of $E$ in the Grothendieck group $K(X)$ of $X$.

\proclaim {Corollary 2.2} Let $ 0 \to E_1 \to E_2 \to E_3 \to 0$ be an exact sequence of coherent sheaves on $X$.  Then $[\lambda^r(E_2)] = \Sigma_{i =0}^r [\lambda^{r-i} (E_1) \otimes_L \lambda^i(E_3)] $ in the Grothendieck group $K(X)$. If the cohomology groups of all the $E_i$ are finite, then $\chi(\tilde \lambda^r(E_2) = \prod_0^r \chi(\tilde \lambda^{r-i}E_1) \chi(\tilde \lambda^iE_3)$ \endproclaim

Proof. Choose  compatible coherent locally free finite resolutions $P^._1$, $P_2^.$, and $P_3^.$ of $E_1$, $E_2$, and $E_3$. Use these resolutions to compute $\tilde \lambda^r$ of $E_1$, $E_2$, and $E_3$ 
Then there exists a filtration on $ \Lambda^rKP_2^.$

$$ 0 =G_0 \subseteq G_1 \subseteq \dots \subseteq G_{r-1} \subseteq G_r = \Lambda^rKP^. _2$$
and exact sequences $ 0 \to G_{i-1} \to G_i \to \Lambda^{r-i}K P^._1 \otimes _s\Lambda^i KP^._3 \to 0 $ .,

 Recall that if if $F^.$ and $G^.$ are simplicial sheaves, the simplicial tensor product of $F^. $ and $G^.$ is given by $(F^. \otimes_s G^.)^n = F^n \otimes G^n$.
 
 We now apply the inverse functor $N$ to our filtration, getting 
 
 $$ 0 = NG_0 \subseteq NG_1 \subseteq \dots \subseteq NG_r-1 \subseteq NG_r = N\Lambda^rKP_2^. = \tilde \lambda^r E_2$$
 and exact sequences $0 \to NG_{i-1} \to NG_i \to N(\Lambda^{r-i}KP_1^. \otimes_s N\Lambda^iKP_3^. )\to 0$. The corollary then follows from the fact ([M]. p. 129ff.  May proves this for simplicial abelian groups, but the argument is valid for any abelian category with tensor products) that if $F^.$ and $G^.$ are simplicial sheaves, $N(F^. \otimes _s G^.)$ is isomorphic in the derived category to $NF^. \otimes _L NG^.$.
 
 \proclaim {Definition 2.3} Let $\lambda^k$ be the usual $k$-th $\lambda$ operation on $K(X)$.  Recall that $\lambda$ is determined by the relations that $\lambda^k([E]) = [\Lambda^k(E)]$ if $E$ is locally free and if $F_2 = F_1 + F_3$ in $K(X)$, then $\lambda^r(F_2) = \Sigma_{i=0}^r \lambda^{r-i}(F_1)\lambda^i(F_3)$.\endproclaim
 
 \proclaim {Theorem 2.4} If $E$ is a coherent sheaf on $X$, $[\tilde \lambda^r(E)] = \lambda^r([E])$. \endproclaim
 
 Proof.  this is an easy double induction on $r$ and the length of a locally free resolution of $E$, using Definition 2.3 and Corollary 2.2

\proclaim {Theorem 2.5} $\chi(\lambda^2(m) )= \chi(k)$, and $\chi(\lambda^2 \kappa(P)) = (\chi(\kappa(P))^{-2}$. \endproclaim

Proof.  Let $P$ be a closed point of $X$, and let $A = O_{X, P}$ be the local ring of $P$ on $X$.  :Let $m$ be the maximal ideal of $A$, and $k$ the residue field of $A$.  Since $A$ is regular, we have the Koszul resolution of $k$ : $0 \to A \to A^2 \to A \to k \to 0$.  

We start with the exact sequence $0 \to m \to A \to k \to 0$.

Since $\lambda^2(A) = 0$ Corollary 2.2 tells us that $\lambda^2(k) + \lambda^2(m) + m \otimes_L k = 0$. 
We have $\chi(m \otimes  _L k) = \chi(A \otimes k) /\chi(k \otimes k)$ which implies because of Theorem 1.3 that $\chi(m \otimes_L k) = \chi (A \otimes_L  k) = \chi (k)$, so $\chi(\lambda^2 k) = (\chi (\lambda^2(m))\chi (k))^{-1}$.

From the exact sequence $0 \to A \to A^2 \to m \to 0$ and Corollary 2.2 we get the triangle in the derived category  $0 \to m \to A \to \lambda^2(m) \to 0$, which implies that $\lambda^2(m) = k$, and hence that $\chi (\lambda^2(k) )= \chi(k) ^{-2}$.

\proclaim {Corollary 2.6} If $F$ is in $Ker(X)$, then $\chi(\lambda^2(F)) = \chi (F)^{-2}$. \endproclaim

Proof.  If we have the exact sequence $0 \to F_1 \to F_2 \to F_3 \to 0$ with $F_i$ finite then $\chi (\lambda ^2(F_2) = \chi(\lambda^2 (F_1) \chi(\lambda^2(F_3)$.  This is an immediate consequence by induction of Corollary 2.2 and Corollary 1.4.  The Corollary then follows immediately from Lemma 1.2.

\S. 3. The conjecture for $d = 1$.

If $d = 1$, $X =$ Spec $O_F$, with $O_F$ the ring of integers in the number field $F$.  The functional equation is well-known and  $A$ is equal to $|d_F|$, where $d_F$ is the discriminant of $F$.  If $ r = o$ our complex $C_{0,1}$ is $O_F \to \Cal D^{-1}$, namely the inclusion of the ring of integers in the inverse different.  This complex has $H^0 = 0$ and $H^1$ isomorphic to $\Omega$, which has order equal to $|d_F|$. So $\chi(C_{0,1} )= A^{-1}$.

If $r = 1$, the complex $C_{1,1}$ is $\Omega \to 0$, so $\chi(C_{1,1}) = A$.

Now let $r \geq 2$. Applying Corollary 2.2  to the sequence $ 0 \to \Cal D \to A \to \Omega \to 0 $, and calling that both $\Lambda^r A$ and $\Lambda^r \Cal D$ are $0$ for $r \geq2$ yields $[\lambda^r \Omega ]+ [\lambda^{r-1}\Omega]= 0$, which of course implies $\chi(\lambda^r [\Omega)] \chi(\lambda^{r-1}[\Omega)] = 1$.

If $r < 0$, the complex $C_{r. 1}$ becomes $0 \to \underline{RHom} (\tilde \lambda^{-r}, \omega)$ (since $\tilde \lambda^r = 0$ if $r <0$) and then the result follows by Serre duality from the cases where $r \geq 1$.

\S. 4 The conjecture for $d =2$

From now on in this paper we will assume $d = 2$.  

\proclaim {Proposition 4.1} Conjecture 0.1 is true for $r =0$ and $r = 1$.\endproclaim

Proof.  If $ r = 1$, the complex $C_{1,2}$ is $\Omega \to \omega$, which is the complex called $C$ in [B].  Bloch proves in [B] that $\chi(C) = A$

If $ r = 0$, $C_{0,2}$ is $ O_X \to \underline{RHom}(\Omega, \omega)$.and the conjecture follows from the conjecture for $r=1$ and Serre duality.

 If $P$ is a closed point of $X$  let $B= O_{X,P}$, $m$ be the maximal ideal of $B$, and $k = \kappa(P)$ be the residue field of $B$.  

\proclaim {Lemma 4.2} Let $r \geq 2$.  Then $\lambda^r[(m)] = (-1)^r [k]$ in $K(X)$, hence  $\chi([\lambda^r(m)]) = (\chi([k]))^{(-1)^r}$. \endproclaim

Proof.  Applying Corollary 2.2 to the exact sequence $0 \to B \to B^2 \to m \to 0$, $r \geq 3$ implies that $[\lambda^r(m) ]+[\lambda^{r-1}(m) ]= 0$, so $\chi[(\lambda^r(m))] \chi[(\lambda^{r-1}(m)) ]= 1$,
and hence since Theorem 2.1 tells us that $\chi([\lambda^2(m)]) = \chi([k])$, Lemma 4.2 follows by induction.

\proclaim {Lemma 4.3} $\chi([\lambda^r(k)]) = (\chi([k]))^r$ if $r$ is odd and $(\chi([k)])^{-r}) $ if $r$ is even.\endproclaim

Proof.  Applying Corollary 2.2 to the exact sequence $0 \to m \to B \to k \to 0$, we obtain $\chi(\lambda^{r+1}(m)) \chi( \lambda^r(k) )\chi( \lambda^{r+1}(k)) = 1$ (Note that Lemma 4.2 implies that $\lambda^j(m)$ is in $Ker (X)$ for $j \geq 2$, and hence Corollary 1.4 implies that $\chi(\lambda^j(m) \otimes_L \lambda^{r-j}(k) )= 1$ for $2 \leq j \leq r-1$.) Induction using Lemma 4.2 now completes the proof.

\proclaim {Proposition 4,4} If $F$ is in $Ker (X)$, then $\chi(\lambda^r(F) = (\chi(F)^r)$ if $r$ is odd and $\chi(F)^{-r}$ if $r$ is even. \endproclaim

Proof, The proof is the same as the proof of Corollary 2.6, starting from Lemma 4.3.

\proclaim{Corollary 4.5} $\chi([\lambda^r([C]))] = A^r$ if $r$ is odd and $A^{-r}$ if $r$ is even. \endproclaim 

Proof.  Since $r_Y(\Omega)$ and $r_Y(\omega$ are both $1$ for all $Y$, $[C]$ is in $Ker(X)$.  Then use Proposition 4.1.

\proclaim{Theorem 4.6} Conjecture 0.1 is true for $r \geq 2$ and $r < 0$.  \endproclaim

Proof.  Since $\omega$ differs from $B$ by something in $Ker(X)$, $\chi([\omega \otimes_L C]) = \chi([B \otimes _L C]) = \chi([C])$.  Then Corollary 2.2 applied to the triangle $C \to \Omega \to \omega \to C[1]$ tells us that if $r \geq 2$, $\chi([\tilde \lambda^r(\Omega) ])= \chi([\tilde \lambda^r(C)]) \chi ([\tilde \lambda^{r-1}(C)])$ and Conjecture 0.1 follows.  The case when $r <0$ follows from Serre duality..  

Appendix (Derived tensor products and derived exterior powers)

The key to defining derived functors (both additive and non-additive) in the absence of projectives is contained in the classic paper [BS] of Borel and Serre.  The basic point  is that if $F$ is a coherent sheaf and  $P^\bullet \to F$ and $Q^\bullet \to F$ are two finite resolutions of $F$ by coherent locally free sheaves, there exists a finite resolution $R^\bullet \to F$ by coherent locally free sheaves which dominates both $P^\bullet \to F$ and $Q^\bullet \to F$. ($R^\bullet \to F$ dominates $P^\bullet \to F$ if there is a surjective map of complexes from $R^\bullet \to F$ to $P^\bullet \to F$ which is the identity on $F$. ) This result follows immediately by induction from Lemma 14 of Borel-Serre.  

Let $L^\bullet $ be the kernel of the map from $R^\bullet \to F$ to $P^\bullet \to F$,  Then $L^\bullet$ is acyclic and locally free , so if $G$ is any coherent sheaf  $L^\bullet \otimes G$ is acyclic.  It immediately follows that the map from $R^\bullet \otimes G$ to $P^\bullet \otimes G$ is a quasi-isomorphism, and so $P^\bullet \otimes G$ and $Q^\bullet \otimes G$ are isomorphic in the derived category.  So we may define the derived tensor product $F\otimes ^L G$ to be $P^\bullet \otimes G$ and this is independent of the choice of the locally free resolution $P^\bullet$.

We wish to define derived exterior powers in an analogous fashion . Let $N$ be the functor which takes simplicial sheaves to bounded below complexes of sheaves, and $K$ be its inverse functor. Let $\Lambda^k$ denote the $k$th exterior power.   We would like to define $\tilde \lambda^k F$ to be $N \Lambda^k K P^\bullet$, where $P^\bullet $  is a resolution of $F$, so we have to show that in the derived category this is independent of the choice of resolution. 

\proclaim {Theorem A,1}  $N\Lambda^k K  P^\bullet$ is independent (in the derived category) of the choice of locally free resolution $P^\bullet$ of $F$. \endproclaim

\proclaim {Lemma A.2}   If $R^\bullet$ is a finite  acyclic complex of locally free sheaves on a noetherian scheme $X$, $\Lambda^k K R^\bullet$ is acyclic for all k $\geq 1$. \endproclaim

Proof . Since exterior powers commute with restriction to open sets, we may assume $X$ is affine.  But then the locally free sheaves are projective, and an acyclic projective complex is homotopically trivial.  The functors $K$, and $\Lambda^k$ preserve homotopy, so the resulting complex is still homotopically trivial, so acyclic. 

Proof of Theorem A.1  Again, by using Lemma 14 of [BS], we reduce to the case where we have a surjective map $f$ from $P^\bullet$ to $Q^\bullet$, where both $P^\bullet$ and $Q^\bullet$ are finite locally free resolutions of the coherent sheaf $F$.  Let $R^\bullet $ be the kernel 0f $f$,  $R^\bullet$ is clearly acyclic, so by the lemma $\Lambda ^i KR^\bullet$ is acyclic, for all $i \geq 1$.   We have the exact sequence of simplicial sheaves $0 \to KR^\bullet \to KP^\bullet  \to KQ^\bullet \to 0 $.  Then $\Lambda^k KP^\bullet$ has a filtration whose associated graded pieces are $\Lambda ^iK R^\bullet \otimes \Lambda^j KQ^\bullet $ for $i + j = k$.  Since $\Lambda ^iKR^\bullet$ is acyclic for $i \geq 1$ we get that $f$ induces a quasi-isomorphism from $\Lambda ^k P^\bullet \to \Lambda^k Q^\bullet$.  Since $N$ preserves quasi-isomorphisms we obtain the desired result.  

References

[B] Bloch, S. De Rham cohomology and conductors of curves. Duke Math. J. 54 (1987), no. 2, 295-"308. 

[BS] Borel, A.; Serre, J.-P. Le thŽorme de Riemann-Roch. Bull. Soc. Math. France 86 (1958) 97-136. 

[H] Hartshorne, R.   Algebraic Geometry,, Springer-Verlag (1977)

[L] Lichtenbaum, S. Special Values of Zeta Functions of Schemes arXiv: 1704.00062

[M] May, Simplicial Objects in Algebraic Topology, Chicago (1967)

[S] Serre, J.-P., Facteurs locaux des fonctions z\`eta des vari\'eti\'es alg\`ebriques (d\'efinitions et conjectures) S\'eminaire Delange-Pisot-Poitou, (1969/70), no. 19 581-592

\bye

\bye